\documentclass{article}

\usepackage{amsfonts}
\usepackage{amssymb}
\usepackage{latexsym}
\usepackage{ifthen}

\newtheorem{thm}{Theorem}
\newtheorem{lem}{Lemma}

\newcommand{\ra}{\rightarrow}
\newcommand{\tx}{\textrm}

\begin{document}

\title{Differential equations defined on algebraic curves}
\author{Vakhtang Lomadze}
\date{~~Andrea Razmadze Mathematical Institute,  Mathematics Department of I. Javakhishvili Tbilisi State University,  Georgia}
\maketitle

{\bf Abstract.}~ The class of ordinary linear constant coefficient differential equations is naturally embedded into a wider class by associating differential equations to   algebraic curves.



\section*{Introduction}

Throughout,  $\mathbb{F}$ is the field of real or complex numbers,  $T$ is an indeterminate, and  $\mathbb{F}\{T\}$ denotes the domain of convergent formal  series in  $T$ with
coefficients in the field $\mathbb{F}$. (We remind that a formal series  $\sum_{n\geq 0}a_nT^n$ is said to be convergent,  if there exists a positive real number  $r$ such that $\sum_{n\geq 0}|a_n|r^n < +\infty$.)
 The fraction field of the domain $\mathbb{F}\{T\}$ is denoted by  $\mathbb{F}(\{T\})$. Elements of this field are referred to as  convergent Laurent series.  By  the theorem on units (see Section 4.4 in Remmert \cite{r}),   stating that an element of the ring  $\mathbb{F}\{T\}$ is invertible
if and only if its free term $\neq 0$ ,
it is immediate  that every convergent Laurent series can be written as a fraction
$f/T^m$, where $f\in \mathbb{F}\{T\}$ and $m\geq 0$.

Throughout, $I$ is an  interval of real axis (with non-empty interior) on which a point is fixed. (The point of this latter  is that it permits
 to define canonically the indefinite integral. In applications, it is customary to regard an interval as a time axis and a fixed point on it as an initial time.)  Without loss of generality, we certainly may assume that $0\in I$  and that  this zero is chosen as a fixed point.
We let $C(I)$ denote the space of all $\mathbb{F}$-valued continuous  functions defined on the interval $I$ and ${\cal D}^\prime(I)$  the space of all $\mathbb{F}$-valued  Schwartz distributions.

The space $C(I)$ has a natural structure of a module over $\mathbb{F}\{T\}$. (The "$T$" acts on continuous functions as the integral operator, and this action is extended in a unique way to all convergent formal series by linearity and continuity.) So that,
for every $f\in \mathbb{F}\{T\}$, we have an operator
$$
C(I)\stackrel{f}{\ra} C(I),\ \ \ \ \ \ w\mapsto fw.
$$
On the other hand, for every $m\geq 0$, there is a differential operator $D^m :{\cal D}^\prime(I)\ra {\cal D}^\prime(I)$. Restricting it to $C(I)$, we get
an operator
$$
 C(I)\stackrel{D^m}{\ra} {\cal D}^\prime(I).
$$
Given a convergent Laurent
series $\phi=f/T^m$ with $f\in \mathbb{F}\{T\}$ and $m\geq 0$, define $D_\phi$ to be the composition
$$
C(I)\stackrel{f}{\ra} C(I)\stackrel{D^m}{\ra} {\cal D}^\prime(I).
$$

A simple observation is  that ordinary linear constant coefficient differential operators have the form $D_\phi$. Indeed,
if $f\in \mathbb{F}[s]$, then the  operator
$$
f(D): C(I)\ra {\cal D}^\prime(I)
$$
can be written as
$$
f(D)=D_\phi,
$$
where $\phi=f(T^{-1})$. Next, as is well-known, one may view $f$ as a rational function on the projective line $\mathbb{P}^1$ associated with
$\mathbb{F}[s]$. Notice also that $\phi$ is the image of $f$ under the embedding of $\mathbb{F}(s)$ (which is the rational function field of $\mathbb{P}^1$)
 into $\mathbb{F}(\{T\})$ that
takes  $s^{-1}$ (which is a local parameter at the infinite point) to $T$.

The theory of linear constant coefficient differential equations can be
 easily generalized to the following situation. We replace $(\mathbb{P}^1,\infty,s^{-1})$ by the triple $(X,P,t)$, where
 $X$ is  an irreducible  smooth projective algebraic curve, $P$ a rational  point on $X$ and $t$ a local parameter of $X$ at this point. The pair $(P,t)$ determines a canonical embedding
$$
\sharp: \mathbb{F}(X)\ra \mathbb{F}(\{T\}),
$$
where $\mathbb{F}(X)$ is the rational function field of $X$.
For every $f\in \mathbb{F}(X)$, we define the operator $f(D): C(I)\ra {\cal D}^\prime(I)$ by setting
$$
f(D)=D_{\sharp(f)}.
$$
It is our belief that differential equations of the form
$$
f(D)w=0
$$
  may be a source for many interesting functions.

To demonstrate that the generalization is meaningful, we show that well-known
Bessel functions are solutions of some differential equations on the hyperbola defined by the equation
 $y^2-x^2=1$. More precisely, we show that, for each nonnegative integer $n$, the  Bessel function
$J_n$ is a solution of the differential equation associated with the rational function $y(x+y)^n$ and that this is the unique solution  that
  satisfies  the following initial conditions
$$
w(0)=w^\prime(0)= \cdots =w^{(n-1)}(0)=0 \ \ \ \tx{and}\ \ \ w^{(n)}(0)=\frac{1}{2^n}.
$$
Moreover, we shall see  that
$$
J_n, D(J_n), \ldots , D^n(J_n)
$$
constitute a fundamental system of solutions.

\section{Differential equations associated with algebraic curves}

For every continuous  function $w\in C(I)$,
let $$\int w$$ denote the "normalized" indefinite integral of $w$ defined by
$$
(\int w)(\xi)=\int_0^\xi w(\alpha)d\alpha,\ \ \ \xi\in I.
$$
Given $f\in \mathbb{F}\{T\}$ and $w\in C(I)$, define the product $fw$ by the formula
 $$
fw=\sum_{n\geq 0} a_n{\int}^n w,
$$
where $a_n$ are the coefficients  of $f$. (The series converges uniformly on every compact neighborhood of $0$.)
This multiplication  makes $C(I)$ a {\em module} over $\mathbb{F}\{T\}$.
 It is worth noting  that
 $$1w=w\ \ \ \ \tx{and}\ \ \ \ Tw =\int w.$$

Define the $E$-transform
$E(f)$ of the convergent formal series $f=\sum_{n\geq 0}a_nT^n$ as the entire analytic function
$$
\xi \mapsto \sum_{n\geq 0}
a_n\xi^n/n!\ \ \ (\xi \in  I).
$$
Remark that
 $
 E(f)=f{\bf 1},
 $
 where ${\bf 1}$ denotes the  constant function that is identically 1 on $I$.

\begin{lem} Let $f,\ g\in \mathbb{F}\{T\}$ and $m,\ n\geq 0$. If \  $T^nf=T^mg$, then the two compositions
$$
C(I) \stackrel{f}{\ra} C(I) \stackrel{D^m}{\ra}  {\cal D}^\prime(I)\ \ \tx{and}\ \ C(I) \stackrel{g}{\ra} C(I) \stackrel{D^n}{\ra}  {\cal D}^\prime(I)
$$
are equal to each other.
\end{lem}
{\em Proof}. This is obvious. Indeed,
$$
\forall w\in C(I),\ \ \ \ \ D^m(fw)=D^{m+n}(T^nfw)=D^{m+n}(T^mgw)=D^n(gw).
$$
$\quad\Box$

If  $f/T^m$ and $g/T^n$ are two  representations of the same convergent Laurent series, then
$T^nf=T^mg$.
In view of the above lemma, it is natural therefore to make the
 following definition.

{\bf Definition}. Let $\phi\in \mathbb{F}(\{T\})$, and assume that $\phi=f/T^m$ with $f\in \mathbb{F}\{T\}$ and $m\geq 0$. Define the  differential operator
$$
D_\phi :  C(I) \ra {\cal D}^\prime(I)
$$
to be the composition
$$
C(I) \stackrel{f}{\ra} C(I) \stackrel{D^m}{\ra}  {\cal D}^\prime(I).
$$

Every $\phi\in \mathbb{F}(\{T\})$  can be uniquely written in the form
$$\phi=\sum_{n>> -\infty}a_nT^n.$$ (The coefficients $a_n$ are zero for all but finitely many negative values of $n$.) The least integer $n$ for which
$a_n\neq 0$ is called the order $\phi$. (In case when $\phi =0$ the order is defined to be $+\infty$.)
Define the degree of the differential equation $$D_\phi w =0$$ as minus  the order of $\phi$. If the degree is non-positive, i.e., if
 $\phi\in \mathbb{F}\{T\}$, then the operator $D_\phi$ is injective, and consequently
the differential equation has no solutions other than 0.
Of interest may be only  equations of positive degree, i.e., those equations that correspond to convergent Laurent  series of the form
$u/T^n$, where $u$   is a unit in $\mathbb{F}\{T\}$  and $n$ is a positive integer.

Let $\mathbb{F}[T]_{\leq k}$ denote the space of polynomials (in $T$) of degree $\leq k$.

\begin{thm} If $\phi\in \mathbb{F}(\{T\})$ and if $\phi=u/T^n$ with  unit $u$ and $n\geq 1$, then
   the  equation
$$
D_\phi w =0
$$
has solutions $$w=E(\frac{p}{u}), \ \ \ p\in \mathbb{F}[T]_{\leq n-1}.$$
\end{thm}
{\em Proof}.  We have:
$$
D_\phi w =0 \ \ \Leftrightarrow\ \  D^n(uw)=0.
$$
 One knows well that the kernel  of the operator
$$
D^n: C(I) \ra  {\cal D}^\prime(I)
$$
is the space of polynomial functions of degree $\leq n-1$, i.e., the space $E(\mathbb{F}[T]_{\leq n-1})$. Because $u$ is invertible, the operator $C(I) \stackrel{u}{\ra} C(I)$ is bijective. Hence,
$$
D^n(uw)=0 \ \Leftrightarrow \ uw\in E(\mathbb{F}[T]_{\leq n-1})\ \Leftrightarrow \ w\in E(u^{-1}\mathbb{F}[T]_{\leq n-1}).
$$
The proof is complete. $\quad\Box$

Of particular interest must be differential equations associated with  those convergent Laurent series that come from algebraic functions, i.e., rational functions on algebraic curves.

Assume we have a triple $(X,P,t)$, where
 $X$ is  an irreducible smooth projective  algebraic curve, $P$ a rational  point on $X$ and $t$ a local parameter of $X$ at this point.
 Letting $O_P$ denote the local ring of $P$ and $\mathbb{F}[[t]]$ the ring of formal series in $t$, we have
 $$
 O_P\subseteq  \mathbb{F}[[t]].
 $$
 In fact, by the implicit function theorem,
 $O_P\subseteq  \mathbb{F}\{t\}$.
 Taking $t$ to $T$, we get a canonical one-to-one ring homomorphism
 $$
 O_P\ra \mathbb{F}\{T\},
 $$
 which, in turn, induces an embedding
 $$
 \mathbb{F}(X)\ra \mathbb{F}(\{T\}),
$$
where $\mathbb{F}(X)$ stands for  the rational function field of $X$.
Denote this canonical embedding by $\sharp$.

{\bf Definition}.
 Given a rational function $f\in  \mathbb{F}(X)$, define a "linear constant coefficient differential" operator $f(D)$ to be $D_{\sharp(f)}$.

We close the section with two examples.

{\em Example 1}. Let $X$ be a projective line, $P$ any its  rational point and $t$ a local parameter at $P$ such that
$s=t^{-1}$ is regular everywhere outside of $P$. Let $\sharp$ be the  embedding of $\mathbb{F}(s)=\mathbb{F}(t)$
into $\mathbb{F}(\{T\})$ for which $\sharp(t)=T$.

An ordinary linear constant coefficient differential equation
$$
a_0w^{(n)} + a_1w^{(n-1)} + \cdots +a_nw = 0
$$
with  $a_0\neq 0$  can be rewritten as
$$
f(D)w=0,
$$
where  $f = a_0s^n + a_1s^{n-1} + \cdots +a_n$. We have
$$
\sharp(f)=T^{-n}(a_0+a_1T+\cdots +a_nT^n).
$$
Consequently, by Theorem 1, the solutions of this equation are
$$
w=E(\frac{p}{a_0+a_1T+\cdots +a_nT^n}),\ \ \ p\in \mathbb{F}[T]_{\leq n-1}.
$$

{\em Example 2}. Let $n$ be a nonnegative integer. Recall that
the Laguerre polynomial of degree $n$ is defined by the formula
$$
L_n(\xi)=\frac{e^\xi}{n!}(\frac{d}{d\xi})^n(\xi^ne^{-\xi}).
$$
Using the Leibniz rule, one has
$$
L_n(\xi)=\frac{e^\xi}{n!}\sum \left(\begin{array}{c}
n\\k\end{array}\right)((-1)^ke^{-\xi})(\frac{n!}{k!}x^k)=\sum \left(\begin{array}{c}
n\\k\end{array}\right)\frac{(-1)^k}{k!}\xi^k,
$$
and thus  $L_n=E((1-T)^n)$.

Let now $(X,P,t)$ and $\sharp$ be as in the previous example. Take the rational function $f=s^{n+1}/(s-1)^n$.
Because
$$
\sharp(s^{n+1}/(s-1)^n)=\frac{T^{-n-1}}{(T^{-1}-1)^n}=T^{-1}(1-T)^{-n},
$$
the solutions of $f(D)w=0$ are
$$
w = cE((1-T)^n)=cL_n,\ \ \ c\in \mathbb{F}.
$$

\section{Bessel functions as special functions associated with hyperbola}

Let us consider the hyperbola defined by the equation
$$
Y^2-X^2=Z^2.
$$
In the affine coordinates $x=X/Z$ and  $y=Y/Z$ this is given by $$y^2-x^2=1.$$
The hyperbola has two points at infinity, namely,
$(1:1:0)$ and $(1:-1:0)$. Take $P=(1:1:0)$, say. In the affine peace $X\neq 0$, the affine coordinates are
$u=Y/X$ and  $t=Z/X$, so that the equation takes the form $u^2-1=t^2$ and our infinite point becomes $(1,0)$.
Choose  $t$ as a local parameter at $P$, and let $\sharp$ denote the canonical embedding
$$
\mathbb{F}(x,y)\ra \mathbb{F}\{T\}
$$
determined by the pair $(P,t)$.
By the implicit function theorem,
$$u=\sqrt{1+t^2}=1 +\frac{1}{2}t^2 + \frac{1}{2}(\frac{1}{2}-1)\frac{t^4}{2!}+\frac{1}{2}(\frac{1}{2}-1)(\frac{1}{2}-2)\frac{t^6}{3!} + \cdots $$
near the point $P$.

Since $y=t^{-1}u$, we have $y=t^{-1}\sqrt{1+t^2}$. For each $n\geq 0$, let us set
$$
j_n=y(x+y)^n.
$$
Letting  $\sqrt{1+T^2}$ denote the convergent formal series
$$
 1 +\frac{1}{2}T^2 + \frac{1}{2}(\frac{1}{2}-1)\frac{T^4}{2!}+\frac{1}{2}(\frac{1}{2}-1)(\frac{1}{2}-2)\frac{T^6}{3!} + \cdots ,
 $$
we have
$$
\sharp(j_n)=T^{-(n+1)}\sqrt{1+T^2}(1+\sqrt{1+T^2})^n.
$$

Remark that $\sqrt{1+T^2}(1+\sqrt{1+T^2})^n$ is a unit in  $\mathbb{F}\{T\}$ and its  free term is equal to $2^n$.
Hence, the  equation $j_n(D)w=0$ has degree $n+1$ and, by Theorem 1, its  solutions  are
$$
w=E(\frac{p}{\sqrt{1+T^2}(1+\sqrt{1+T^2})^n}),\ \ \ p\in \mathbb{F}[T]_{\leq n}.
$$

Put
$$
J_n= E(\frac{T^n}{\sqrt{1+T^2}(1+\sqrt{1+T^2})^n}).
$$
One sees that this is the solution with initial conditions
$$
w(0)=w^\prime(0)= \ldots w^{(n-1)}(0)=0 \ \ \ \tx{and}\ \ \ w^{(n)}(0)=\frac{1}{2^n}.
$$
One can see also that
$$
J_n, D(J_n), \ldots , D^n(J_n)
$$
constitute a fundamental system of solutions.

We are going to show that the functions $J_n$ are none other than the Bessel functions.

\begin{lem} We have:
$$
J_0(\xi)=1-\frac{\xi^2}{4}+\frac{\xi^4}{(2!)^24^2}-\frac{\xi^6}{(3!)^24^3}+ \cdots
$$
and
$$
J_1(\xi)=\frac{\xi}{2}-\frac{\xi^3}{2\cdot 2^3}+\frac{\xi^5}{2!\cdot 3!\cdot 2^5}-\frac{\xi^7}{3!\cdot 4!\cdot 2^7}+ \cdots.
$$
\end{lem}
{\em Proof}.
Using the equality
$$1\cdot 3\cdot \ldots \cdot (2k-1)=\frac{(2k)!}{k!2^k},$$
we find that
$$
\frac{1}{\sqrt{1+T^2}}=\sum_{k\geq 0} \frac{(-\frac{1}{2})(-\frac{1}{2}-1)\ldots (-\frac{1}{2}-k+1)}{k!}T^{2k}=  \ \ \ \ \ \ \ \ \ \ \ \ \ \ \ \ \ \ \
$$
$$\ \ \ \ \ \ \ \ \ \ \ \ \ \ \ \ \ \ \ \
\sum_{k\geq 0} \frac{(-1)^k\cdot 1\cdot 3\cdot \ldots \cdot (2k-1)}{2^kk!}t^{2k}=\sum_{k\geq 0} \frac{(-1)^k (2k)!}{4^k(k!)^2}T^{2k}.
$$
It follows  that
$$
J_0(\xi)=E(\frac{1}{\sqrt{1+T^2}})= \sum_{k\geq 0} \frac{(-1)^k}{4^k(k!)^2}\xi^{2k}.
$$

Next, we have
$$
\frac{T}{\sqrt{1+T^2}(1+\sqrt{1+T^2})}=\frac{\sqrt{1+T^2}-1}{T\sqrt{1+T^2}} =T^{-1}(1-\frac{1}{\sqrt{1+T^2}})= \ \ \ \ \ \ \ \ \ \ \ \ \ \ \
$$
$$
\sum_{k\geq 1} \frac{(-1)^{k+1} 1\cdot 3\cdot \ldots \cdot (2k-1)}{2^kk!}T^{2k-1}=\sum_{k\geq 1} \frac{(-1)^{k+1} (2k)!}{2^{2k}k!k!}T^{2k-1}=
\sum_{k\geq 1} \frac{(-1)^{k+1} (2k-1)!}{2^{2k-1}(k-1)!k!}T^{2k-1}.
$$
It follows that
$$
J_1(\xi)=E(\frac{T}{\sqrt{1+T^2}(1+\sqrt{1+T^2})})=\sum_{k\geq 1} \frac{(-1)^{k+1}}{2^{2k-1}(k-1)!k!}\xi^{2k-1}.
$$
The proof is complete. $\quad\Box$

\begin{lem} For every  $n\geq 1$, we have:
$$
2J^\prime_n=J_{n-1}-J_{n+1}.
$$
\end{lem}
{\em Proof}. For $n\geq 0$, put $$g_n=\frac{T^n}{\sqrt{1+T^2}(1+\sqrt{1+T^2})^n}.$$
We have
$$
2(1+\sqrt{1+T^2})=(1+\sqrt{1+T^2})^2-T^2.
$$
Multiplying this by $T^{-1}g_{n+1}$, we get
$$
2T^{-1}g_n=g_{n-1}-g_{n+1}.
$$
Applying the $E$-transform, we complete the proof. $\quad\Box$

It immediately follows from the above two lemmas that $J_n$ are Bessel functions (of first kind).

Thus, we have proved the following theorem.

\begin{thm} Let $n$ be a nonnegative integer. A fundamental system of solutions  of the differential equation
$$
j_n(D)w=0
$$
is
$$
J_n, D(J_n), \ldots , D^n(J_n),
$$
where $J_n$ is the Bessel function of order $n$. Moreover, $J_n$ is the solution
that satisfies the following initial conditions
$$
w(0)=w^\prime(0)= \cdots =w^{(n-1)}(0)=0 \ \ \ \tx{and}\ \ \ w^{(n)}(0)=\frac{1}{2^n}.
$$
\end{thm}

\end{document}